\documentstyle[12pt]{article}
\newcommand{\gppres}[2]{\left<#1\left|#2\right.\right>}
\newcommand{\ra}{\rightarrow}
\newcommand{\setcond}[2]{\left\{#1\left|#2\right.\right\}}
\newcommand{\csastargp}{CSA$^\ast$-group}
\newtheorem{theorem}{Theorem}
\newtheorem{prop}{Proposition}
\newtheorem{lemma}{Lemma}
\newtheorem{coro}{Corollary}
\newcommand{\gstar}{G^\ast}
\newtheorem{defn}{Definition}
\newcommand{\subgp}[1]{\left<#1\right>}
\newcommand{\length}[1]{\left|#1\right|}
\newcommand{\ben}{\begin{enumerate}}
\newcommand{\een}{\end{enumerate}}
\newcommand{\beq}{\begin{equation}}
\newcommand{\eeq}{\end{equation}}
\newcommand{\boldz}{\mbox{{\bf Z}}}
\newcommand{\boldq}{\mbox{{\bf Q}}}
\newcommand{\boldr}{\mbox{{\bf R}}}
\newcommand{\te}[1]{t^{e_{#1}}}

\newcommand{\bea}{\begin{eqnarray*}}
\newcommand{\eea}{\end{eqnarray*}}
\newcommand{\set}[1]{\left\{\,#1\,\right\}}
\newcommand{\calb}{{\cal B}}
\newcommand{\sdp}{\mbox{$>\hspace{-.16cm}\lhd$}}
\newcommand{\lra}{\longrightarrow}
\newcommand{\csastar}{CSA$^\ast$}
\newcommand{\ultprod}{\prod_{\Phi}^\ast}
\begin{document}
\bibliographystyle{c:/pctex/texbib/plain}
\title{CSA groups and separated free constructions}
\author{D. Gildenhuys, O. Kharlampovich and A. Myasnikov}
\maketitle

\section{Introduction}

The class of hyperbolic groups and the class of groups that can act freely
on a $\Lambda$-tree are two rapidly developing areas of group theory that
have attracted the attention of specialists from many fields of mathematics.
Their study involves an interplay of geometry and algebra. Another very
active area of
group theory, this time connected to mathematical logic, is the study of
systems of equations in free groups and hyperbolic groups. The success of
this
study has inspired renewed interest in Tarski's famous problem: ``is the
elementary theory of a free nonabelian group algorithmically decidable?''
One of the key steps in this direction is the description of $\exists$-free
groups (i.e., those groups that have the same existential theory as a
non-abelian free group). The purpose of this article is to study the class
of CSA-groups, which contains the 
classes of torsion-free hyperbolic groups \cite{Grom},
groups acting freely on $\Lambda$-trees \cite{BassNAT} 
and $\exists$-free groups \cite{vnr}, \cite{GagSpell}.
CSA-groups share many of the properties of the groups in the above-mentioned
three classes, but have the advantage of being definable in purely
group-theoretic terms.

\begin{defn} \cite{BBNewman} We define a subgroup $H$ of a group $G$ to be 
{\em malnormal} (also called {\em conjugate separated})
if $H\cap H^x=1$ for all $x\in G-H$.\end{defn}

It is clear that the intersection of a family of malnormal subgroups is again
malnormal, which allows us to define the {\em malnormal closure}
$\mbox{mal}_G(A)$ of a
subgroup $A$ of a group $G$ to be the intersection of all the malnormal
subgroups of $G$ containing $A$.

\begin{defn}\cite{MyasExpo2}
A group $G$ is called a {\em CSA-group} if all its
maximal abelian subgroups are conjugate separated.\end{defn}

\begin{defn}
 Let $G$ be a group and $\phi : A \rightarrow B $
an isomorphism of subgroups of $G$. The HNN-extension
$$ \gstar = \gppres{G,t}{t^{-1}at =\phi(a),a\in A} $$
is called:
\begin{itemize}
\item[1)] {\em separated} if $A\cap g^{-1}Bg =1$ for all $g\in G;$
\item[2)] {\em strictly separated} if $A\cap g^{-1}\mbox{mal}_G(B)g =1$
for all $g\in G$.
\end{itemize}

\label{DefSepHNN} \end{defn}

O. Kharlampovich and A. Myasnikov recently proved \cite{KhMyasHyp} that the
class
of hyperbolic groups is closed under separated HNN-extensions, subject to
the condition that the two associated subgroups $A$ and $B$ be
quasi-isometrically embedded in $G$ and one of them be malnormal. 
The class of groups acting freely on
$\Lambda$-trees is also closed under separated HNN-extensions, provided that
$A$ and $B$ are abelian and satisfy some natural compatibility condition
\cite{MyasRemLen}.

In the same spirit, we prove in section~\ref{HNN} 
that separated HNN-extensions of an arbitrary
\csastar-group (CSA, without elements of order two) with associated
malnormal subgroups  is again \csastar. In fact, a much more general
result (Theorem~\ref{SepExt}) is true: any strictly separated HNN-extension
of a \csastar-group $G$ 
with associated subgroups $A$ and $B$ is again \csastar\ if  $A$ is
malnormal in $G$ and $B$ is normal in $\mbox{mal}_G(B)$.
We obtain similar results
for amalgamated products of \csastar-groups
(Theorem~\ref{amalgprod}) and, in Theorem~\ref{GraphGroups}, for the
fundamental groups of certain types of graphs of \csastar-groups.
We show that in Theorem~\ref{SepExt}, which is the  main result of the 
section, neither the requirement that the HNN-extension be separable,
nor the requirement that
one of the associated subgroups be malnormal, is a
necessary condition for the HNN-extension to be CSA
(this is not astonishing, because even non-abelian free
group can be decomposed in very complicated ways in terms of HNN-extensions and
free products with amalgamation). However, in the most important case, 
where the associated subgroups are abelian, we
use our main result to obtain in the next section a complete description
of those HNN-extensions which preserve the \csastar-property.

We start section~\ref{AssAb} by showing that
if an HNN-extension with abelian associated subgroups is CSA, then
at least one of the associated subgroups must be maximal abelian. 
Under these assumptions there then
remain four types of HNN-extensions, among which only the separated
HNN-extensions and the rank~1 extensions of centralizers
(described in Proposition~\ref{ConjExt})
preserve the \csastar-property.
Theorem~\ref{AbSepExt}, which states the preservation of the \csastar-property
by separated HNN-extensions, has as its corollary a similar result about
amalgamated products and tree products of \csastar-groups.

Gersten has conjectured that a torsion-free one-relator group is
hyperbolic if and only if it does not contain any Baumslag-Solitar groups
$$B_{m,n}=\gppres{x,y}{yx^my^{-1}=x^n},\ mn\neq 0 $$
(HNN-extensions of the infinite
cyclic group). We prove that a torsion-free one-relator group
fails to be
CSA if and only if it contains a nonabelian metabelian
Baumslag-Solitar group $B_{1,n}$, $n \neq 1$, or the group 
$\calb=F_2\times \boldz$,
the product  of a free group on two generators by
an infinite cycle (Theorem~\ref{OneRelCSA});
a one-relator group with torsion fails to be a CSA-group
if and only if it contains the infinite dihedral group $D_\infty$
(Theorem~\ref{OneRelT}). This gives a complete description, in terms of
``obstacles'', of one-relator CSA-groups.
Every one-relator group with torsion is hyperbolic. This has already
been observed in \cite{BGSS}; it follows
easily from a theorem of B. B. Newman (see \cite{BBNewman} or
\cite{LynSchupp}, Proposition~5.28 on p.~109).

It is known that CSA-groups are commutative transitive (a group $G$ is
commutative transitive if the relation ``$a$ commutes with $b$'' is
transitive on the set $G-\set{1}$;
i.e., if the centralizer of every nontrivial
element of $G$ is abelian, \cite{MyasExpo2}, Proposition~10). That the
converse is not true was pointed out in \cite{MyasExpo2}, and is shown by
the simple example of the infinite dihedral group $D_\infty$.
However, it {\em is} true for torsion-free one-relator groups
(Theorem~\ref{CommTrans}). It is not true,
in general, for one-relator groups with torsion,
because all these groups are commutative transitive
(B.B.Newman \cite{BBNewman}, Theorem~2) but some of them are not CSA
(see the example in Proposition~\ref{TObstacles}).

In section~\ref{ExpGps} we consider
exponential groups and the property of being residually of prime power order
(``residually $p$''). Let $A$ be an (associative) ring with
identity. A group $G$ is called an
{\em $A$-group} if its elements admit exponents from the ring $A$. The
defining axioms can be found in \cite{MyasNova}, \cite{MyasExpo1},
\cite{MyasExpo2}. The {\em tensor completion over $A$} of a group is defined
by the obvious universal property, and a group is said to be {\em
$A$-faithful} if the canonical morphism from the group into its tensor
completion over $A$ is injective. Tensor completions of groups have been
studied extensively in \cite{bau}, \cite{MyasExpo1} and \cite{MyasExpo2}.
Myasnikov and
Remeslennikov proved that if $G$ is a torsion-free CSA-group and $A$ a ring
whose underlying abelian group is torsion free, then $G$ is $A$-faithful,
and the tensor completion of $G$ over $A$ is again a torsion-free CSA-group
(\cite{MyasExpo2}, Theorem~9). So, the class of torsion-free CSA-groups is
contained in the class of $\boldq$-faithful groups
(it is clear that a $\boldq$-faithful group must be torsion free).
We prove
(Proposition~\ref{aapfaith}) that if, for almost all primes $p$, a group
is residually $p$, then it is $\boldq$-faithful.

Gilbert Baumslag \cite{Ba2} posed the problem
of determining which one-relator groups are $\boldq$-faithful.
The class of torsion-free one-relator CSA-groups is strictly contained in
the class of one-relator $\boldq$-faithful groups as well as in the class
of one-relator groups that are residually $p$ for almost all primes $p$
(Proposition~\ref{OneRelNotCSA}).
It would be interesting to find obstacles
to $\boldq$-faithfulness in torsion-free one-relator groups (in the sense
that the groups $\cal B$ and $B_{1,n}$ of Theorem~\ref{OneRelCSA} are
obstacles to the CSA-property, and the groups $B_{m,n}$ are conjectured to
be obstacles to hyperbolicity).

In Proposition~\ref{OneRelNotCSA}, we give an example of a one-relator group
which is not CSA, but is nevertheless $\boldq$-faithful, residually
torsion-free nilpotent and residually $p$, with torsion-free
pro-$p$-completion, for every prime $p$.
 On the other hand metabelian non-abelian Baumslag-Solitar 
groups provide examples of ${\bf Q}$-faithful
groups which are non-CSA and not residually $p$ for almost all primes $p$.

\section{HNN-extensions of \csastargp s}\label{HNN}

Our main purpose in this section is to establish some natural sufficient
conditions for HNN-extensions and amalgamated products to preserve
the CSA-property. We also
investigate the more general problem of determining when the fundamental group
of a graph of CSA-groups is again CSA.

\begin{defn}A {\em \csastargp} is a CSA-group without
elements of order 2.\end{defn}

Throughout the paper, $\gstar$ will denote an HNN-extension
of a group $G$ relative to an isomorphism of associated subgroups
$\phi: A \rightarrow B$, and $t$ will denote the stable letter.
Recall that any HNN-extension $\gstar$ of a group $G$
is endowed with a length function (\cite{LynSchupp}, p. 185).
We denote the length of an element $z$ of $\gstar$ by $\left|z\right|$.

\begin{lemma} Let $G^*$ be a strictly separated HNN-extension of a group $G$
with associated subgroups $A$ and $B$ such that $A=\mbox{mal}_G(A)$ and
$B \unlhd \mbox{mal}_G(B)=B_1$. Let $c\in\gstar$.\ben
\item If $1\neq b_1\in B_1$, $b_1^c\in G$ and $c\notin G$, then
$c\in B_1t^{-1}G$ and $b_1\in B$.\label{ConjInBOne}
\item If $B_1^c\cap B_1\neq 1$, then\label{BOneInt}
$c\in B_1$.
\item If $1\neq a\in A$, $a^c\in G$, then there are three possibilities:
$c\in A$, $c\in tG$ or $c\in tB_1t^{-1}G$.\label{ConjInA}
\item If $A^c\cap A\neq 1$, then $c\in tB_1t^{-1}$.\label{AInt}
\een\label{NiceLemma}  \end{lemma}
\underline{Proof}. Suppose that $c\notin G$. We write $c$ in reduced form:
\beq c = g_0\te{1}g_1\te{2}g_2\cdots g_{n-1}\te{n}g_n,\ \left|c\right|=n\geq 1.
\label{redformc} \eeq
\ref{ConjInBOne}. Suppose that $b_1\in B_1$ and $b_1^c\in G$. We have
\beq c^{-1}b_1c = g_n^{-1}t^{-e_n}g_{n-1}^{-1}\cdots
g_2^{-1}t^{-e_2}g_1^{-1}t^{-e_1}g_0^{-1}b_1g_0\te{1}g_1\te{2}g_2\cdots
g_{n-1}\te{n}g_n.\label{conjbone}\eeq
Since the length of the left side of this equation is 0, the right side is
not reduced and
$$ g_0^{-1}b_1g_0 \in \left\{\begin{array}{l}
A\mbox{ if } e_1=1 \\
B\mbox{ if } e_1=-1 \end{array}\right.$$
If $e_1=1$, then $b_1\in A^{g_0^{-1}}\cap B_1$, which violates our
hypotheses. Hence, $e_1=-1$, $g_0^{-1}b_1g_0\in B$,
$b_1\in B_1\cap B^{g_0^{-1}}\subseteq B_1\cap B_1^{g_0^{-1}}$, which implies
that $g_0\in B_1$ and $b_1\in B^{g_0^{-1}}\subseteq B$. Write
$a=tb_1^{g_0}t^{-1}$.
If $n\geq 2$, then we see from (\ref{conjbone}) that
the word $t^{-e_2}g_1^{-1}ag_1\te{2}$ is not reduced, and
$$ g_1^{-1}ag_1 \in \left\{\begin{array}{l}
B\mbox{ if } e_2=-1 \\
A\mbox{ if } e_2=1.\end{array}\right.$$
If $e_2=-1$ then $a\in B\cap A^{g_1^{-1}}\subseteq B_1\cap A^{g_1^{-1}}$,
which contradicts our hypotheses. So, $e_2=1$,
$a\in A\cap A^{g_1^{-1}}$ and $g_1\in A$. But then this contradicts the
assumption that (\ref{redformc}) was in reduced form. It follows
that $n=1$ and $c=g_0t^{-1}g_1\in B_1t^{-1}G$.\\
\ref{BOneInt}. Since $B_1$ is malnormal in $G$,
it is sufficient to derive a contradiction from our assumption that
$c\notin G$. It follows from part~1 that if $1\neq b_1\in B_1\cap cB_1c^{-1}$,
then $b_1\in B$ and for some $g\in G$, $c\in B_1t^{-1}g$. But then
$$ b_1^c \in g^{-1}tBt^{-1}g\cap B_1=A^g\cap B_1 ,$$
which contradicts our assumptions.\\
\ref{ConjInA}. We keep our assumption that $c$ has reduced form
(\ref{redformc}) (if $c\in G$ and $a^c\in A$, then $c\in A$, since $A$ is
malnormal). By a similar reasoning as in part~\ref{ConjInBOne} we find that
$e_1=1$, $g_0\in A$. Letting $a^{g_0t}$ and $t^{-1}g_0^{-1}c$ play the roles
of $b_1$ and $c$ respectively in part~\ref{ConjInBOne}, we find that
$t^{-1}g_0^{-1}c\in G$ or $t^{-1}g_0^{-1}c\in B_1t^{-1}G$; i.e.,
$c\in g_0tG=tG$ or $c\in g_0tB_1t^{-1}G=tB_1t^{-1}G$.\\
\ref{AInt}. If $A^c\cap A\neq 1$, then $B_1^{c^t}\cap B_1\neq 1$ and
$c\in tB_1t^{-1}$ by part~\ref{BOneInt}.\\
This completes the proof of the lemma.

\begin{lemma} Let $G^*$ be a strictly separated HNN-extension of a group $G$
with associated subgroups $A$ and $B$ such that $A=\mbox{mal}_G(A)$ and
$B \unlhd \mbox{mal}_G(B)=B_1$. Suppose that $M$ is a maximal abelian subgroup 
of
$\gstar$. Then one of the following is true.\ben
\item The intersection of $M$ with some conjugate $A^{x^{-1}}$ ($x\in\gstar$)
of $A$ is nontrivial, in which case $M\subseteq xtB_1t^{-1}x^{-1}$.
\item The intersection of $M$ with every conjugate of $A$ is trivial, but
$M$ intersects some conjugate of $G$ nontrivially, in which case
$$ \left(\forall x\in\gstar\right)\left(M\cap G^x\neq 1
\Rightarrow M\subseteq G^x\right). $$
\item There exists a cyclically reduced element $z\in\gstar$ such that
$M$ is a conjugate of the infinite cyclic group $\subgp{z}$,
$\length{z}\geq 1$.\een

\label{maxabingstar}\end{lemma}

\underline{Proof}.
It follows from \cite{SerreTrees}, Chapter 1, section 5.4,
Theorem~13 and its Corollary~2 that $\gstar$ can be made to act on a tree
$X$ in such a way that the stabilizers of the edges are conjugates of $A$
in $\gstar$ (the vertices of $X$ are cosets $yG$ of $G$ in $\gstar$,
and the edges
are cosets of $A$ in $\gstar$). If a nontrivial element of $M$
stabilizes some edge of $X$, then there exists and element $x\in\gstar$ such
that $M\cap A^{x^{-1}}\neq 1$, and hence $M^{xt}\cap B_1\neq 1$. Every element
of $M^{xt}$ then centralizes a nontrivial element of $B_1$, and it follows
from Lemma~\ref{NiceLemma}, part 2, that $M^{xt}$ is contained in $B_1$,
and hence $M$ in $xtB_1t^{-1}x^{-1}$.

Suppose now that $M$ intersects trivially all the stabilizers of edges of
the tree $X$. It
follows from \cite{SerreTrees}, Chapter 1, section 5.4 Theorem 13
(see also Example 1) of section 5.1) that $M$ is a free product of
conjugates of subgroups of $G$ and a free group. Since $M$ is abelian, the
free product is trivial (i.e., it has only one factor). More precisely,
if a nontrivial element of
$M$ fixes some vertex of $X$, say the coset wG (equivalently
$M\cap wGw^{-1}\neq 1$), then $M$ is entirely
contained in the conjugate $wGw^{-1}$ of $G$; i.e., $M$ is a conjugate of a
maximal abelian subgroup of $G$.

On the other hand, if $M$ acts
freely on the tree $X$ then $M$ is free, hence
cyclic, and is generated
by an element $w$ of length $\geq 1$. Since $w$ is not conjugate to an
element of $G$ ($w$ does not fix any vertex of $X$), it is conjugate to a
cyclically reduced element $z$, with $\left|z\right|\geq 1$. This completes
the proof of the lemma.

\begin{theorem} Let $G^*$ be a strictly separated HNN-extension of 
a \csastar\ group $G$ with associated subgroups $A$ and $B$ such that 
$A=\mbox{mal}_G(A)$ and $B \unlhd \mbox{mal}_G(B)$. Then $\gstar$ is a 
\csastargp.
\label{SepExt}\end{theorem}

\underline{Proof}. It is clear that $\gstar$ has no elements of order~2
(\cite{LynSchupp}, Theorem~2.4, p.~185).
We may assume that $A$ and $B$ are nontrivial (else
$\gstar$ is the free product of $G$ and an infinite cyclic, hence CSA
\cite{MyasExpo2}).
As usual we will use the notation $B_1=\mbox{mal}_G(B)$.

Let $M$ be a maximal abelian subgroup of $\gstar$.
Suppose that $\exists v\in\gstar$ such that
$M\cap M^v\neq 1$. We must prove that $v\in M$.

We first consider the case where some element of $M$
has nontrivial intersection with a conjugate $A^{x^{-1}}$ of $A$
($x\in\gstar$). Then,
according to Lemma~\ref{maxabingstar}, $M^{xt}\subseteq B_1$.
The groups $M^{xt}$ and $M^{vxt}$ have a nontrivial element $w$ in common,
which belongs to $B_1$ and which is centralized by $M^{vxt}$. By
Lemma~\ref{NiceLemma}, $M^{vxt}$ (which can be written $M^{xtv^{xt}}$)
is  contained in $B_1$ and
$$ 1\neq w\in B_1\cap B_1^{v^{xt}}.$$
It follows from Lemma~\ref{NiceLemma} that $v^{xt}\in B_1\subseteq G$.
Since $G$ is CSA, and
$$ 1\neq w\in M^{xt}\cap M^{xtv^{xt}} ,$$
we conclude that $v^{xt}\in M^{xt}$ and $v\in M$.

We now assume that $M$ intersects trivially every
conjugate of $A$ but intersects nontrivially a conjugate
$G^w$ of $G$, $w\in \gstar$. In this case, Lemma~\ref{maxabingstar}
tells us that $M$ is entirely contained in $G^w$. But then $M^v$ also
intersects $G^w$ nontrivially, and there exists a maximal abelian subgroup
$N$ of $G$ such that $M=M^v=N^w$ ($G$ is commutative transitive). Next, we
let $v'=wvw^{-1}$, then $N^{v'}=N$, and we
claim that $v'\in G$. If not, we write it in reduced form:
\beq v' = v_0\te{1}v_1\te{2}v_2\cdots v_{n-1}\te{n}v_n,\ 
\left|v'\right|=n\geq 1.
\label{redformv} \eeq
Let $z\in N$. We have
\beq \left(v'\right)^{-1}zv' = v_n^{-1}t^{-e_n}v_{n-1}^{-1}\cdots
v_2^{-1}t^{-e_2}v_1^{-1}t^{-e_1}v_0^{-1}zv_0\te{1}v_1\te{2}v_2\cdots
v_{n-1}\te{n}v_n.\label{Conjz}\eeq
Since the length of the left side of this equation is 0, the right side is
not reduced and
$$ v_0^{-1}zv_0 \in \left\{\begin{array}{l}
A\mbox{ if } e_1=1 \\
B\mbox{ if } e_1=-1 \end{array}\right.$$
In both cases we find that $M$ has a nontrivial intersection with a conjugate
of $A$ in $\gstar$, contrary to our assumption that it intersects trivially
all the stabilizers of edges of the tree $X$. So $v\in G$.
Since $G$ is a CSA-group, $v\in M$.

Suppose now that the third possibility
of Lemma~\ref{maxabingstar} applies to $M$. Replacing $M$ by one of its
conjugates, we may suppose that $M=\subgp{z}$, with $z$ cyclically reduced.
All powers of $z$ are then cyclically reduced
as well. There exist integers $m,n$ such that
$v^{-1}z^mv=z^n$. By the Conjugacy Theorem for HNN-extensions
(\cite{LynSchupp}, Chapter~4, Th. 2.5),
$$ |m|\length{z} = \length{z^m} = \length{z^n} = |n|\length{z}; $$
hence $m=\pm n$ and $v^2$ lies in the center of the group $\subgp{z^n,v}$.
Suppose that $\subgp{z^n,v}$ intersects a conjugate $A^{x^{-1}}$ of $A$
nontrivially, then $\exists w\in \subgp{v,z^n}^{xt}\cap B_1$, $w\neq 1$.
Since $\left(v^{xt}\right)^2$ centralizes $w$,
it follows from Lemma~\ref{NiceLemma} that
$\left(v^{xt}\right)^2\in B_1$. Since $\left(z^{xt}\right)^n$ centralizes the
nontrivial element $\left(v^{xt}\right)^2$ of $B_1$,
it follows again from Lemma~\ref{NiceLemma} that
$\left(z^{xt}\right)^n\in B_1$. But $z^n$, being cyclically reduced and of
length $\geq 1$, cannot belong to a conjugate of $G$. This proves that
$\subgp{z^n,v}$ is a free product of subgroups of conjugates of $G$ and a
free group. Since $\subgp{z^n,v}$ has nontrivial center, it is indecompasable
as a free product (\cite{MagKarSol}, Section 4.1, Problem~21, p.195). Since,
as already pointed out, $z^n$ is not conjugate to an element of $G$,
$\subgp{z^n,v}$ is contained in a free group, hence must be cyclic. This
implies that $z^n$ is in the center of $\subgp{z,v}$. By the same argument
as before, if $\subgp{z,v}$ intersects some conjugate of $A$ then $z$ is
contained in a conjugate of $G$, which is impossible, since $z$ is cyclically
reduced of length $\geq 1$. So $\subgp{z,v}$ is a free product of subgroups
of conjugates of $G$ and a free group. It has nontrivial center, which implies
as before that $\subgp{z,v}$ is indecomposable, and hence contained in a free
group. This implies that $v$ commutes with $z$. By the maximality of
$M=\subgp{z}$, $v\in M$. This is what we had to prove.

\begin{coro} A separated  HNN-extension of a \csastargp\  with malnormal
associated subgroups is a \csastargp.
\end{coro} 

Similar results for amalgamated products can be easily obtained using
the following lemma.

\begin{lemma} Let $A$ and $B$ be subgroups of groups $G$ and $H$ respectively,
$\phi : A \rightarrow B $ an isomorphism. The groups $A$ and $B$ can be 
considered as isomorphic subgroups of the free product $G \ast H$. Let us 
denote by
$$E\left(G,H,\phi\right) = \gppres{G \ast H, t}{ t^{-1}at=\phi(a)}$$
the HNN-extension, associated with $\phi$, of the group $G \ast H$.
The amalgamated product $G\ast_\phi H$ is embeddable in
$E\left(G,H,\phi\right)$.\label{AmalFromHNN}
\end{lemma}

\underline{Proof}. The subgroup $\subgp{G^t,H}$ generated in
$E\left(G,H,\phi\right)$ by the
$t$-conjugate of $G$ and $H$ is isomorphic to $G\ast_\phi H$. It can be 
easily verified using the normal forms of elements
in $E\left(G,H,\phi\right)$.

\begin{theorem} Let $G$ and $H$ be \csastargp s,
$A$ and $B$  subgroups
of $G$ and $H$, respectively, such that $A=\mbox{mal}_G(A)$ and
$B \unlhd \mbox{mal}_H(B)$,
and $\phi:A\ra B$ an isomorphism. Then the
amalgamated product $G\ast_\phi H$ is \csastar.\label{amalgprod}
\end{theorem}

\underline{Proof}.
Let $B_1=\mbox{mal}_H(B)$. First, we claim that $A$ and $B_1$ are also
malnormal
in the free product $G\ast H$. Malnormality is transitive: if $X$ is a
malnormal subgroup of $Y$ and $Y$ is a malnormal subgroup of $Z$, then $X$
is a malnormal subgroup of $Z$. So we need only point out (as Lyndon and
Schupp already did on page~203 of \cite{LynSchupp}) that the factors of a free
product are malnormal in the product. Then we need the fact that $A$ and
$B$ are mutually conjugate separated, i.e. $A\cap B_1^x=1$ for
all $x\in G\ast H$. It is enough to prove that $G$ and $H$ are mutually
conjugate separated in their free product, and this is easily verified using
normal forms. Next, we need the fact that the class
of \csastargp s is closed under free products (\cite{MyasExpo2}, Theorem~4).
To complete the proof, we note that
the amalgamated product $G\ast_\phi H$ is embedded in the HNN-extension of
$G\ast H$, relative to the isomorphism $\phi$ (Lemma~\ref{AmalFromHNN}),
and we apply
Theorem~\ref{SepExt} (note that subgroups of \csastargp s are \csastar,
\cite{MyasExpo2}, Proposition~13).

\begin{coro} An amalgamated product of \csastargp s with malnormal amalgamated
 subgroups is again \csastargp.
\end{coro} 

To deal with graphs of groups we need the next Proposition,
mentioned without proof in \cite{BBNewman}.

\begin{prop} Let $A$ be a malnormal subgroup of a group $G$, $B$ a subgroup
of a group $H$, $\phi:A\ra B$ an isomorphism, and
$P=G\ast_\phi H$ the associated amalgamated
product. Then every malnormal subgroup of $H$ is also
malnormal in $P$.\label{MalInAmal}\end{prop}

\underline{Proof}. To simplify the exposition, we will suppose that $G$ and
$H$ are subgroups of $P$ and $\phi$ is the identity.
Because of the transitivity of malnormality, it suffices
to prove that $H$ is malnormal in $P$.
Every nontrivial element $x$ of $P$ can be written in the form
$$ x = p_1p_2\cdots p_r, $$
where each $p_i$ lies in one of the factors $G$ or $H$, and no $p_i$ lies
in $A$ if $r>1$; also, if the length $r$
of $x$ is $>1$ then $p_i$ and $p_{i+1}$ lie in different factors
($i=1,\ldots,r-1$). Although this representation for $x$ is not necessarily
unique, the length $r$ is unique, and so is the sequence of factors
determined by $p_1,p_2,\ldots,p_r$.
Suppose that $h\in H\cap H^x$. We shall prove by induction on $r$ that
$p_i\in H$ for all $i=1,\ldots,r$. There exists an $h'\in H$ such that
$$ h = p_r^{-1}\cdots p_2^{-1}p_1^{-1}h'p_1p_2\cdots p_r.$$
It follows that $p_1^{-1}h'p_1\in A$ or $r=1$. \\
\underline{Case 1}: $r=1$. If $p_1\in H$, then $x\in H$ and we are done.
If $p_1\in G$, then $p_1^{-1}h'p_1h^{-1}=1$. The left side is reducible,
which implies that $h'\in A$ and $h\in H\cap G=A$. But, $A$ is malnormal
in $G$, so $x=p_1\in A$. So $p_1\in H$.\\
\underline{Case 2}: $r>1$. Let $a_1=p_1^{-1}h'p_1\in A$. Then $p_1\in H$
by Case~1. \\
We can now write
$$ h = p_r^{-1}\cdots p_2^{-1}a_1p_2\cdots p_r,$$
and apply the induction hypothesis to complete the proof.

Note that under the hypotheses of Proposition~\ref{MalInAmal} the malnormal
closure of $A$ in $P$ contains the normalizer of $B$ in $H$.

Following Dicks \cite{DicksGG} we define an {\em oriented graph of groups}
as follows. It consists of an oriented graph
$\Gamma=\left(V,E,\bar{\iota},\bar{\tau}\right)$ ($\bar{\iota}(e)\in V$ and
$\bar{\tau}(e)\in V$ are the initial and terminal vertices respectively of
an edge $e\in E$), together with a function $G$
which assigns to each vertex $v\in V$ a group $G(v)$, and to each edge
$e\in E$ a subgroup $G(e)$ of $G\left(\bar{\iota}(e)\right)$ and monomorphism
$t_e:G(e)\ra G\left(\bar{\tau}(e)\right)$. We shall say that this oriented
graph of groups is 
{\em quasi-malnormal} if, for all $e\in E$, $G(e)$ is malnormal in
$G\left(\bar{\iota}(e)\right)$ while $t_e(G(e))$ is normal in its malnormal
closure in $G\left(\bar{\tau}(e)\right)$.  If, in addition, for each edge
$e$, $t_e(G(e))$ is malnormal in $G\left(\bar{\tau}(e)\right)$, then we say
that the graphs of groups is {\em malnormal}. In this case, the
orientation of the graph is irrelevant. The graph of groups is said to be 
{\em separated} if for any edge $e$ which is
a loop ($v=\bar{\iota}(e)=\bar{\tau}(e)$) one has 
$G(e)^g\cap t_e\left(G(e)\right)=1$  for all $g\in G(v)$.

If a quasi-malnormal separated oriented graph of \csastar-groups has only 
one edge, then its fundamental group is \csastar\ (Theorem~\ref{SepExt},
Theorem~\ref{amalgprod}, \cite{MyasExpo2}, Theorem~6).
We prove in
Proposition~\ref{BadTree} that the fundamental group of a quasi-malnormal 
oriented tree of CSA-groups need not be CSA.
For every ordinal number $\alpha$ we define as follows an oriented graph
$L_\alpha$ (a ``line'' from 1 to $\alpha$). Its vertices are ordinals
$\leq\alpha$, and $\setcond{\left(\beta,\beta +1\right)}
{1\leq \beta,\beta+1\leq\alpha}$ is its set
of edges. The edge $\left(\beta,\beta +1\right)$ has initial vertex $\beta$,
and terminal vertex $\beta+1$.

\begin{theorem} The fundamental group of an oriented
graphs of \csastar-groups is again \csastar\ in the following
two situations:\ben
\item the underlying oriented graph is $L_\alpha$ for some ordinal $\alpha$,
and the oriented graph of groups is quasi-malnormal;\label{line}
\item the underlying oriented graph is a tree and the graph of groups is
malnormal.\label{tree}\een\label{GraphGroups}\end{theorem}

\underline{Proof}.
\ref{line}. The class of \csastar-groups is closed under direct limits
(\cite{MyasExpo2}, Theorem~6) and free products (Theorem~\ref{amalgprod}),
hence we are reduced to proving the result
for finite $\alpha$. The result then follows, by a simple induction argument,
from Theorem~\ref{amalgprod} and Proposition~\ref{MalInAmal}.\\
\ref{tree}. By a similar argument, we need only prove the result for finite
trees, and the result follows from the Theorem and Proposition.

The separation conditions, necessary to ensure that the fundamental group
of an arbitrary malnormal separated oriented graph of \csastargp s is
\csastar, are quite complicated and cumbersome to formulate.

\underline{Remark}. The conditions of Theorem~\ref{SepExt} are not
necessary: 
\begin{itemize}
\item[1)] the following example $\gstar$ of an HNN-extension of
a \csastar-group $G$ is \csastar, without
the HNN-extension being separated (or
a centralizer extension, in the sense of Definition~6, \cite{MyasExpo2}).
$$ \gstar = \gppres{x_1,x_2,x_3,t}
{x_1^t=x_2, x_2^t=x_1x_3},$$
$G$ is the free group on $x_1,x_2,x_3$, and the associated subgroups
$A=\subgp{x_1,x_2}$, $B=\subgp{x_2,x_1x_3}$, are malnormal in $G$,
but have nontrivial intersection $\subgp{x_2}$. Clearly, $\gstar$ is
free on two generators $x_1, t$, hence CSA.
\item [2)] Let
$$ \gstar = \gppres{x_1,x_2,x_3,x_4,x_5,t}{x_1^t=x_1^{x_2}, 
\left(x_1^{x_2}\right)^t=x_3, x_4^t=x_5^2}.$$
Then this is a representation of $\gstar$ as a non-separated HNN-extension
of the free group $G=\subgp{x_1, x_2, x_3, x_4, x_5}$,
with associated subgroups $A=\subgp{x_1, x_1^{x_2}, x_4}$ and
$B=\subgp{x_1^{x_2}, x_3, x_5^2}$ (both $A$ and $B$ are free of rank~3).
The subgroups $A$ and $B$ are not malnormal in $G$; moreover, neither is
normal in its malnormal closure. The group $\gstar$ can be
represented as a free product with amalgamation:
$$\gstar=\gppres{x_1,x_2,x_3,x_4,t}{x_1^t=x_1^{x_2},
\left(x_1^{x_2}\right)^t=x_3} \ast_{\psi}\subgp{x_5}, $$
where $\psi: \subgp{x_4^t} \rightarrow \subgp{x_5^2}$ is the obvious
isomorphism of infinite cycles. Let us denote the left factor of the
above decomposition by $G_1$.
The subgroup $\subgp{x_4^t}$ is maximal abelian in $G_1$, therefore, to prove
that $\gstar$ is a CSA-group
it is enough to prove that $G_1$ is a CSA-group (Theorem~\ref{amalgprod}).
Writing $s=tx_2^{-1}$, one can represent $G_1$ as follows:
$$ G_1 = \gppres{x_1,x_2,x_4,s}{\left[x_1,s\right]=1}.$$
Hence the group $G_1$ is CSA (see Proposition~\ref{ConjExt}).
\end{itemize}

It is interesting that both examples above can be constructed 
using only 
``admissable'' HNN-extensions and free products with amalgamation (i.e. those
mentioned in our theorems as sufficient conditions).

\section{HNN extensions of \csastargp s with abelian associated subgroups}
\label{AssAb}

In this section we give a complete description of HNN-extensions, with abelian
associated subgroups, that preserve the \csastar-property.

\underline{Remark}. An abelian subgroup of a CSA-group is malnormal iff it 
is maximal abelian.

\begin{prop} Let $\gstar$ be an HNN-extension of a CSA-group $G$, relative
to an isomorphism $\phi:A\ra B$ of nontrivial abelian subgroups $A$ and $B$
of $G$. Then the HNN-extension $\gstar$ is strictly separated if and only
if it is separated. \label{AbSep}\end{prop}

\underline{Proof}.
Suppose that $A^s\cap B_1\neq 1$, $s\in G$; then the CSA-property of $G$
implies that the maximal subgroups $A^s$ and $B_1$ of $G$ are equal. Hence,
$A^s\cap B\neq 1$. The result follows.

\begin{prop} \ben\item Suppose that $\phi:A\ra B$
is an isomorphism of two nontrivial
abelian subgroups $A$ and $B$ of a group $G$. If neither $A$ nor $B$ is
maximal abelian in $G$ then the associated HNN-extension $\gstar$ is not
CSA.
\item Suppose that $\phi:A\ra B$
is an isomorphism of two nontrivial abelian subgroups of two groups $G$ and
$H$ respectively. If $A$ is not maximal abelian in $G$ and $B$ is not maximal
abelian in $H$, then the associated amalgamated product is not CSA.\een
\label{MustMax}\end{prop}

\underline{Proof}. 1. Suppose that $\gstar$ is CSA.
Let $A_1$ and $B_1$ be maximal abelian subgroups of $G$
containing $A$ and $B$ respectively. Let $a_1\in A_1-A$, $b_1\in B_1-B$ and
let $t$ denote as usual the stable letter of the HNN-extension $\gstar$,
so that $a^t=\phi(a)$ for all $a\in A$. Let $1\neq b\in B$. Then $a_1^t$
commutes with $b$, and so does $b_1$. Since CSA-groups are commutation
transitive,
$$ 1=\left[a_1^t,b_1\right] = t^{-1}a_1^{-1}tb_1^{-1}t^{-1}a_1tb_1 .$$
This is impossible, since the word on the right cannot be reduced.\\
2. A similar argument can be used in this case.

Suppose as before that $\phi:A\ra B$ is an isomorphism of two nontrivial
abelian subgroups $A$
and $B$ of a CSA-group $G$, with $A$ maximal abelian in $G$. Then one of the
following four possibilities must apply:\begin{itemize}
\item[1)] $\left(\forall s\in G\right)\left(A^s\cap B=1\right)$;
\item[2)] $\left(\exists s\in G\right)\left(A^s=B\,\wedge\,
\left(\forall a\in A\right)\left(\phi(a)=a^s\right)\right)$;
\item[3)] $B$ is maximal abelian in $G$,
$(\exists v\in G)\left(A^v\cap B\neq 1\right)$
and
$$ \left(\forall s\in G\right)\left[A^s\cap B\neq 1\Rightarrow
A^s=B\,\wedge\, \left(\exists a_0\in A\right)\left(a_0^s\neq
\phi\left(a_0\right)\right)\right];$$
\item[4)] $B$ is not maximal abelian in $G$, and $\left(\exists s\in G\right)
\left(A^s\cap B\neq 1\right)$ (in this case, $B\subset A^s$).
\end{itemize}
In the HNN-extension
$$ \gstar = \gppres{G,t}{t^{-1}at =\phi(a),a\in A}, $$
$\phi(a)=a^t$ for all $a\in A$. In cases 3) and 4), there exist $s\in G$
and $a_0\in A$ such that
$$ A^{ts^{-1}}\subseteq A \mbox{ but } \left[a_0,ts^{-1}\right]\neq 1,\
\mbox{and hence }ts^{-1}\notin A. $$
This shows that in the cases 3) and 4) the group $\gstar$ is not CSA.

\begin{prop} Let $A$ be a maximal abelian subgroup of a \csastargp\
$G$, $B$ a subgroup of $G$ and $\phi:A\ra B$ an isomorphism. Suppose that
there exists an $s\in G$ such that $B=A^s$ and $\phi(a)=a^s$ for all
$a\in A$. Then the HNN-extension
$\gstar=\gppres{G,t}{t^{-1}at=\phi(a)}$ is a \csastargp.
\label{ConjExt}\end{prop}

\underline{Proof}. It is clear that if we put $v=ts^{-1}$, then
$$ \gstar =\gppres{G,v}{\left[a,v\right]=1,\forall a\in A}.$$
I.e., we get a presentation for $\gstar$ from a presentation for $G$
by taking as generators for $\gstar$ those of $G$ together with $v$, and
as defining relations those of $G$ together with $\left[a,v\right]=1$ for
all $a\in A$. Since $A$ is maximal abelian, we have $A=C_G(a)$ for some fixed
$a\in A$, and $\gstar$ is the direct, rank 1, extension of the centralizer
of the element $a$ (Definition 6, \cite{MyasExpo2}). Since $A\times\subgp{v}$
has no element of order 2, it follows from Theorem 5 \cite{MyasExpo2} that
$\gstar$ is a \csastargp.

\begin{theorem}
Let $G$ be a \csastargp\ and $\gstar$ a separated HNN extension of $G$,
relative to an isomorphism $\phi:A\ra B$ of abelian subgroups of $G$, with
$A$ maximal abelian in $G$.
Then $\gstar$ is a \csastargp.\label{AbSepExt}\end{theorem}

\underline{Proof}. We may assume that $A$ and $B$ are nontrivial (else
$\gstar$ is the free product of $G$ and an infinite cyclic, hence \csastar\
\cite{MyasExpo2}). The maximal abelian subgroup $B_1$ of $G$ containing $B$
is the malnormal closure of $B$. Clearly, $B$ is normal in $B_1$,
and the result follows from Theorem~\ref{SepExt} and Proposition~\ref{AbSep}.

Combining Proposition~\ref{ConjExt}, Theorem~\ref{AbSepExt} and the remarks
at the beginning of this section, we obtain

\begin{theorem} Let $\phi:A\ra B$ be an isomorphism of abelian subgroups
of a \csastar-group $G$, with $A$ maximal abelian in $G$. Then the corresponding
HNN-extension of $G$ is again \csastar\ if and
only if it is a separated extension
or there exists an $s\in G$ such that $B=A^s$ and $\phi(a)=a^s$ for all
$a\in A$. \end{theorem}

\begin{theorem} Let $G$ and $H$ be \csastargp s and
$\phi:A\ra B$ an isomorphism of abelian subgroups of $G$ and $H$ respectively.
Then the amalgamated product $G\ast_\phi H$ is \csastar\ if and only if at
least one of the subgroups $A$ or $B$ is maximal abelian in $G$ or $H$
respectively.\label{amalgiff}
\end{theorem}

\underline{Proof}. It follows directly from Proposition~\ref{MustMax} and
Theorem~\ref{amalgprod}.

\underline{Remark}. An oriented graph $\Gamma$ of CSA-groups,
with abelian edge groups,
is quasi-malnormal if and only if each edge group $G(e)$ is maximal abelian in
$G\left(\bar{\iota}(e)\right)$; and $\Gamma$ is malnormal if and
only if each edge group $G(e)$ is maximal abelian in both vertex groups
(i.e., $G(e)$ is maximal abelian in
$G\left(\bar{\iota}(e)\right)$ and its image $t_e\left(G(e)\right)$ is maximal
abelian in $G\left(\bar{\tau}(e)\right)$). Theorem~\ref{GraphGroups} gives
examples of types of oriented graphs of \csastar-groups, with abelian edge
groups, whose fundamental groups are again \csastar. We have, in particular,

\begin{prop} If the edge groups $G(e)$ of an oriented tree of \csastar-groups
are maximal abelian in $G\left(\bar{\iota}(e)\right)$ and also have maximal
abelian images in the target groups $G\left(\bar{\tau}(e)\right)$, then
the fundamental group of the tree of groups is again a \csastar-group.
\label{TreeProdAb}\end{prop}

\underline{Proof}. The result follows immediately from
Theorem~\ref{GraphGroups} and the
Remark at the beginning of this section.

\begin{prop} Let $T=\left(V,E\right)$ be the tree with
two edges $e_1$ and $e_2$, having a common initial point $v$, and endpoints
$w_1$ and $w_2$ respectively. Suppose that the vertex groups are \csastar\ and
the edge groups $G\left(e_i\right)$ are maximal abelian in $G(v)$, with
nontrivial intersection (hence coincident). If their images are not maximal
abelian in $G\left(w_1\right)$ and $G\left(w_2\right)$ respectively, then
the fundamental group of this oriented graph of groups is not CSA.
\label{BadTree}\end{prop}

\underline{Proof}. The fundamental group is the amalgamated product of the
CSA-groups
$G(v)\ast_{G\left(e_1\right)} G\left(w_1\right)$
(Theorem~\ref{amalgprod}) and $G\left(w_2\right)$,
with amalgamated subgroup $G\left(e_2\right)$. The maximal abelian subgroup
of $G_v\ast_{G\left(e_1\right)} G\left(w_1\right)$, containing the image
of $G\left(e_2\right)$, contains the image of
the maximal abelian subgroup of
$G\left(w_1\right)$ containing the image of $G\left(e_1\right)$, hence cannot
coincide with the image of $G\left(e_2\right)$. The result now follows from
Proposition~\ref{MustMax}.

\section{One-relator CSA-groups}
\label{OneRel}

In this section we completely characterise, in terms of ``obstacles'', all
one-relator groups that are CSA, and we show that all the obstacles are
realized. It follows from the characterization of
one-relator CSA-groups that in the torsion-free case the CSA-property is
equivalent to the transitivity of commutation; however, this equivalence
fails for one-relator groups with torsion.

\begin{prop} \ben
\item The group $\calb = F_2\times\boldz$ is not commutative transitive
(hence not CSA).   \label{calb}
\item The non-abelian Baumslag-Solitar groups\label{nonabBS}
$$B_{m,n}=\gppres{x,y}{yx^my^{-1}=x^n},\  mn\neq 1,$$
are not CSA; 
furthermore, the non-metabelian Baumslag-Solitar groups
($\left|m\right|\neq 1\neq \left|n\right|$) contain
$\calb$ as a subgroup (hence are not commutative transitive and not CSA)
\item The one-relator group
$$ G = \gppres{x,y}{\left[\left[x,y\right],y\right]=1} $$
contains $\calb$, but does not contain any non-abelian Baumslag-Solitar
groups.\label{TripleComm}
\een\label{TFObstacles}\end{prop}

\underline{Proof}. \ref{calb}. If we write
$\calb=F\left(x,y\right)\times\subgp{z}$, then we see immediately that the
elements $x$ and $y$ commute with $z$, but not with each other.\\
\ref{nonabBS}.
A metabelian non-abelian BS-group contains a nontrivial abelian normal 
subgroup, therefore it is not a CSA-group.
See the proof of Theorem~\ref{OneRelCSA} for a verification
of the statement that the non-metabelian Baumslag-Solitar groups contain
$\calb$. \\
\ref{TripleComm}. The proof that $\calb$ is a subgroup of $G$ is contained in
the proof of Proposition~\ref{OneRelNotCSA} of the next section.
Since $G$ is residually $p$ (Proposition~\ref{OneRelNotCSA}) for all primes
$p$ but the non-abelian Baumslag-Solitar groups fail to be residually $p$
for almost all $p$ (Proposition~\ref{res}), $G$ does not
contain any non-abelian Buamslag-Solitar groups.

\begin{lemma} \label{sdprods} Suppose that a group $T$ is a nonabelian
semidirect product
$A^+\sdp \subgp{t}$ of the underlying abelian group $A^+$ of a finitely 
generated
subring $A$ of $\boldq$ and an infinite cycle $\subgp{t}$. Then every
epimorphism $\psi$ from $T$ onto a torsion-free group $W$, with $\psi$
nontrivial on $A$ and $\subgp{t}$, is an isomorphism.\end{lemma}

\underline{Proof}. It is easily seen that $A=\boldz\left[\frac{1}{k}\right]$
for some natural number $k$, and that the restriction of $\psi$ to $A^+$
is injective. The action of $t$ on $A^+$ is multiplication by some
rational number $\frac{m}{n}$, where both $m$ and $n$ divide a power of $k$.
Clearly,
$\psi(T)=\psi\left(A^+\right)\left<\psi(t)\right>$. Identify
$\psi\left(A^+\right)$ with $A^+$. The action (through
conjugation) by $\bar{t}$ on $A^+$ is multiplication by $\frac{m}{n}$. If
no power of $\bar{t}$ is in $A^+$, then $\psi$ is clearly an isomorphism.
Suppose now that some power $\bar{t}^i$ is in $A^+$. Then, since $A^+$ is
abelian, $m^i=n^i$ and $m=\pm n$. Since $T$ is nonabelian, $\frac{m}{n}=-1$.
The $\subgp{t}$-module $A^+$ is the union of an ascending chain of submodules
$\subgp{a_i}$; therefore, $T$ is the union of an ascending chain of groups,
each isomorphic to $H=\gppres{a,t}{t^{-1}at=a^{-1}}$. If the restriction of
$\psi$ to each of these groups is injective then $\psi$ is injective; so
we may assume, without
loss of generality, that $T$ is in fact equal to $H$. If $\psi$ is
not injective then $\psi(T)$ admits a presentation with generators $a,t$,
and
$$ t^{-1}at = a^{-1},\ t^r = a^s\ (r,s\in\boldz)$$
among its defining relations. Then
$$ t^r = t^{-1}a^st = \left(t^{-1}at\right)^s = a^{-s} = t^{-r} $$
and $t^{2r}=1$. This contradicts our assumption that $\psi(T)$ is torsion
free.

\begin{theorem}\label{OneRelCSA}
A nonabelian torsion-free one-relator group $G$ is
CSA if and only if it does not
contain a copy of $\calb$ or one of the (nonabelian) metabelian
Baumslag-Solitar groups $B_{1,n}=\gppres{x,y}{yxy^{-1}=x^n}$,
$n\in\boldz-\set{0,1}$.
\end{theorem}

The proof of Theorem~\ref{OneRelCSA} will use the following

\begin{lemma}
Let $H$ be a subgroup of a torsion-free one-relator group $G$. Then one of
the following is true:\begin{itemize}
\item $H$ is locally cyclic;
\item $H$ contains a nonabelian free group of rank~2;
\item $H$ is isomorphic to $B_{1,m}$ for some $m\in\boldz-\set{0}$.
\end{itemize}
\label{subgplemma}\end{lemma}

\underline{Proof of the lemma}.
By \cite{LynSchupp}, Chapter II, Proposition~5.27, $H$ is either solvable
or contains a free group of rank~2. Suppose that $H$ is solvable. The
result then follows from Moldavanskii's Theorem
(\cite{LynSchupp}, Chapter II, Prop. 5.25, p. 109) or from the
classification of solvable groups of cohomological dimension $\leq 2$
\cite{DGSol1} ($G$ has cohomological dimension $\leq 2$, and the same must
be true of its subgroups).

\underline{Proof of the Theorem}. Clearly $G$ cannot be CSA if it contains
a copy of $\calb$ or $B_{1,n}$.
If $G$ is not
CSA, then there exists a maximal abelian subgroup $A$ of $G$, and elements
$a_1,a_2\in A$, $z\in G-A$, such that $a_1^z=a_2\neq a_1$.
By Lemma~\ref{subgplemma}, $\subgp{a_1,a_2}$ is either free abelian of
rank two or it is cyclic.

Consider first the case where $\subgp{a_1,a_2}$ is cyclic. Then
$\exists a_0\in A$, $m,n\in\boldz$ such that $a_1=a_0^m$, $a_2=a_0^n$.
If $\left|m\right|=1$, then, by Lemma~\ref{sdprods},
$\subgp{a_0,z}$ is isomorphic to the semidirect
product of the additive group of the ring $\boldz\left[\frac{1}{n}\right]$
and the infinite cycle $\subgp{z}$, where the action by $z$ is multiplication
by $n$. Similarly, if $\left|n\right|=1$ then $\subgp{a_0,z}$ is isomorphic
to the semidirect product of $\boldz\left[\frac{1}{m}\right]$ and the infinite
cycle $\subgp{z}$, where the action of $z$ is now multiplication by $m$.
So, if $\length{n}$ or $\length{m}$ is 1, then $\subgp{a_0,z}$
is isomorphic to $B_{1,n}$ or $B_{1,m}$. It cannot be isomorphic to the free
abelian group of rank~2, $B_{1,1}$, since we have supposed that
$a_1\neq a_2$. Suppose now that
$\left|m\right|\neq 1\neq \left|n\right|$.
Consider, for every natural number $i$, the group
$$ H_i=\subgp{a_0^{z^{-i}},a_0^{z^{-i+1}},\ldots,
a_0^{z^{-1}},a_0,a_0^{z},\ldots,a_0^{z^i}}.$$
These groups cannot all be abelian, since in that case $\subgp{a_0,z}$ would be
isomorphic to a homomorphic image of the semidirect product of
$\boldz\left[\frac{1}{mn}\right]$
and $\subgp{z}$, where the action of $z$ is multiplication by $\frac{m}{n}$;
hence, by Lemma~\ref{sdprods}, $\subgp{a_0,z}$ would be isomorphic to this
(metabelian) semidirect product. But
this solvable group is not among those that appear in the statement of
Lemma~\ref{subgplemma} (it has cohomological dimension 3, \cite{DGSol1}).
It follows from the Lemma that either some $H_i$ is isomorphic to $B_{1,m}$,
$m\notin\set{0,1}$, or $H_i$ contains a free group of rank~2.
It is easily proved, by induction on $i\geq 1$, that
\beq a_0^{(mn)^i} = \left(a_0^{m^{2i}}\right)^{z^i}
= \left(a_0^{n^{2i}}\right)^{z^{-i}}. \label{PowerConj}\eeq
So, $a_0^{(mn)^i}$ centralizes $H_i$, and it follows that, if $H^i$ contains
a free group of rank~2, then
$\subgp{a_0,z}$ contains a copy of the group $\calb$. In particular, the
above arguments show that
$B_{m,n}$ contains a copy of $\calb$.

We now consider the case where
$\subgp{a_1,a_2}$ is not cyclic. Suppose that $G$ does not contain any
nonabelian Baumslag-Solitar group. Then every subgroup of $G$ is either
abelian or contains a free group of rank~2. Let $K_i$ be the subgroup of
$G$ generated by $a_1^{z^j}$, $-i\leq j\leq i$. If all $K_i$ are abelian,
then their union $K$ is abelian and normal in the nonabelian solvable subgroup
$\subgp{K,z}$ of $G$. By the Lemma and our assumptions, $G$ cannot contain
a nonabelian solvable subgroup. Hence, we may assume that $K_i$
is abelian and $K_{i+1}$ is not. Then $K_{i+1}$ contains a free subgroup
of rank~2, and we claim that $K_0$ centralizes it. If $i=0$, this is clear
since $\subgp{a_1,a_1^z}$ and $\subgp{a_1^{z^{-1}},a_1}$ are abelian. If
$i>0$, then
$$ \subgp{a_1,a_1^z,\ldots,a_1^{z^{i+1}} },\
\subgp{a_1^{z^{-i-1}},\ldots,a_1^{z^{-1}},a_1} $$
are contained in conjugates of $K_i$, hence are abelian.
The result follows.

\begin{prop}
The infinite dihedral group
$$D_\infty=\gppres{x,y}{yxy^{-1}=x^{-1}, y^2=1}$$
is commutative transitive, but not CSA. Moreover, $D_\infty$ is a subgroup
of the one-relator group
\beq G = \gppres{x,y}{x^2=1} .\label{CTNotCSA}\eeq
\label{TObstacles}\end{prop}

\underline{Proof}. Clearly,
$\subgp{x}\cap\subgp{x}^y=\subgp{x}$, $y\notin\subgp{x}$, which shows
that $D_\infty$ cannot be CSA.

The normal subgroup of $G$ generated by $x$ is the free product of countably
many copies of a cyclic group of order 2, hence contains $D_\infty$.

We can represent $G$ as a free product of $\boldz$ and $\boldz/2\boldz$.
To prove that $G$, and hence $D_\infty$, is commutative transitive, it
suffices to observe that a free product $P$ of commutative transitive groups
is commutative transitive. Indeed, if $x$ commutes with $y$ and $z$ in $P$,
then $H=\subgp{x,y,z}$ is a free product of conjugates of subgroups of the
factors of $P$ and (possibly) a free group $F$. Since the center of
$H$ is nontrivial, the decomposition of $H$ as a free product is trivial,
hence $H$ is isomorphic to a subgroup of $F$ or of one of the factors of $P$.
Thus, $y$ commutes with $z$.

\begin{theorem}\label{OneRelT}
Let $G$ be a one-relator group with torsion. Then $G$ is CSA if and only
if it does not contain the infinite dihedral group
$$ D_\infty = \gppres{x,y}{yxy^{-1}=x^{-1}, y^2=1}.  $$
\end{theorem}

\underline{Proof}. If $G$ contains $D_\infty$, it cannot be CSA, since the
class of CSA-groups is closed under taking subgroups and $D_\infty$ is not
CSA, by Proposition~\ref{TObstacles}.

Conversely, suppose that $G$ is not CSA and does not contain a copy of
$D_\infty$. By a result of Karass and Solitar, every subgroup of a
one-relator group with torsion is either cyclic, $D_\infty$ or contains a
free group of rank~2 (\cite{LynSchupp}, Chapter II, Prop.~5.27).
Let $A$ be a maximal abelian subgroup of $G$. Then $A$ is cyclic.

\underline{Case 1}: $A$ is infinite.\\
Suppose
that $1\neq a_1\in A=\subgp{a_0}$, $a_1=a_0^m$, $z\in G$, $a_2=a_1^z=a_0^n$,
$z\notin A$, then there are two possibilities for
$$ K=\gppres{a_0^{z^i}}{i\in\boldz}:$$
it contains a free subgroup of rank~2 or is cyclic.
If $K$ is cyclic, then $\subgp{K,z}$ does not contain a nonabelian free
group. Hence is abelian. It contains the maximal abelian subgroup $A$ of
$G$, hence $z\in A$. Contradiction.
If $K$ contains a free subgroup $F$
of rank~2, then $F$ is contained in the subgroup generated by a
finite number of the displayed generators of $K$; hence is centralized by
a power $a_0^{(mn)^i}$ of $a_0$ (see the proof of Theorem~\ref{OneRelCSA}).
Since the center of $F$ is trivial, any nontrivial element of $F$, together
with $a_0^{(mn)^i}$, generate a noncylic abelian subgroup, and we again have a
contradiction.\\
\underline{Case 2}: $A$ is finite.\\
We apply the usual Magnus treatment to the one-relator
group (\cite{LynSchupp}, Chapter~2, section~6, Chapter~4, section~5), and
argue by induction on the length of the defining relator. If only one letter
appears in the defining relator, then the one-relator group is finite cyclic,
or the free product of a free group and a finite cyclic group. It is easily
verified that in this case the one-relator group is again CSA (the base
of the induction). If the exponentsum of no generator in the
relator is zero, we can adjoin a root of a generator to our one-relator group
$G$, and take a different system of generators, so that the relator, when
expressed in these new generators, has exponentsum
zero in one of the generators, and is an HNN-extension of a one-relator group
$H$, with shorter defining relator. The associated subgroups, the co-called
Magnus subgroups, are free, hence have trivial intersection with $A$.
It follows that $A$ is a free product of subgroups of conjugates of $H$ and
a free group. But $A$, being finite cyclic, is indecomposable as a free
product, and must be contained in a conjugate of $H$. As in the
proof of Theorem~\ref{SepExt} (see the discussion around (\ref{redformv})
and (\ref{Conjz})), we find that $z$ belongs to the same conjugate of $H$.
By the induction hypothesis, $H$ is CSA. But then $z\in A$, the desired
contradiction.

\begin{theorem} A torsion-free one-relator group is CSA if and only if it
is commutative transitive. The class of one-relator groups with torsion that
are CSA is strictly contained in the class of commutation-transitive
one-relator groups with torsion.\label{CommTrans}\end{theorem}

\underline{Proof}. In the torsion-free case, the obstacles $\calb$ and
$B_{1,n}$ ($n\neq 1$) to the CSA-property are not commutative transitive.
The result follows. In the torsion case the group $G$ of
Proposition~\ref{TObstacles} is commutative transitive but not CSA.

\section{\boldq-faithfulness, residual properties and one-relator groups}
\label{ExpGps}

In the context of G. Baumslag's problem \cite{Ba2} of describing the class
of $\boldq$-faithful one-relator groups, we prove here that this class
strictly contains the union of the class of one-relator \csastargp s
and the class of one-relator groups that are residually $p$ for almost all
primes $p$.

\begin{prop} If, for almost all primes $p$, a group $G$ is residually $p$
then it is $\boldq$-faithful.\label{aapfaith}\end{prop}

\underline{Proof}. Let $S$ be the set of prime numbers and $\Phi$ an
ultrafilter on $S$, containing all cofinite subsets of $S$. The ultraproduct
$\ultprod\boldz_p$ of the rings of $p$-adic numbers contains $\boldq$ as
a subring (every integer is divisible by only finitely many primes).
Moreover, the componentwise action by exponentiation of this ultraproduct on
the ultraproduct $\ultprod \hat{G}_p$ of pro-$p$-completions $\hat{G}_p$
of $G$ is faithful. The restriction to $\boldq$ of this action is faithful.
Since $G$ is embeddable in the $\boldq$-group $\ultprod \hat{G}_p$, $G$ is
$\boldq$-faithful.

\begin{prop} If a group $G$ is residually torsion-free nilpotent, then it
is faithful over $\boldq$.\label{tfnilp}\end{prop}

\underline{Proof}.
By hypothesis, $G$ is
embeddable in a product of torsion-free nilpotent groups. Since torsion-free
nilpotent groups are faithful over $\boldq$ (see \cite{MalNilp} or
\cite{HallNilp}), each one is embeddable in its
$\boldq$-completion, and it follows that there is a monomorphism from $G$
into a product of $\boldq$-groups. The canonical homomorphism from $G$
into its $\boldq$-completion factors through this monomorphism, hence is
itself a monomorphism, which means that $G$ is $\boldq$-faithful.

\begin{prop} Every nonabelian Baumslag-Solitar group $B_{m,n}$
is not residually $p$ for almost all primes $p$. If
$B_{m,n}$ is metabelian (i.e. $\left|m\right|$ or $\left|n\right|=1$), then
it is $\boldq$-faithful.\label{res}
\end{prop}

\underline{Proof}. The defining relation of a nonabelian Baumslag-Solitar
$G=B_{m,n}$ can be written in the form
\beq x^{n-m} = \left[x^m,y^{-1}\right] .\label{CommForm}\eeq
Consider first the case were $n\neq m$.
If $p$ is a prime not dividing $m-n$, then $G$ is not residually
$p$. Indeed, the pro-$p$-completion $\hat{G}_p$ is a one-relator pro-$p$-group
with the same defining relation, and if $x$ belongs to the $k$-th term of
the central descending series of $\hat{G}_p$, then (\ref{CommForm}) shows that
$x^{n-m}$, and hence $x$, belongs to the $(k+1)$-st term. This proves that
the image of $x$ in $\hat{G}_p$ is 1, and the canonical map from $G$ into
$\hat{G}_p$ cannot be injective, which means that $G$ is not residually $p$.
If $m=n$ then $y$ commutes with $x^n$ but not with $x$. However, for every
prime $p$ not dividing $n$ the pro-$p$-completion $\hat{G}_p$ of $G$ is an
abelian pro-$p$-group. This shows that the canonical map:$G\ra \hat{G}_p$
is not injective, and $G$ is not residually $p$.

Clearly, the nonabelian metabelian groups are semidirect products
of the form $\boldz\left[\frac{1}{n}\right]\sdp \boldz$, where the action
of a generator of $\boldz$ on $\boldz\left[\frac{1}{n}\right]$ is multiplication 
by $n>1$.
This semidirect product is naturally embedded in the semidirect product
$E=\boldr\sdp\boldq$, where the action
$$ \theta:\boldq\lra\mbox{Aut}\,\boldr $$
is given by $\theta\left(\frac{s}{t}\right)(r)=n^{\frac{s}{t}}r$. For every
$\left(c,d\right)\in E$ and natural number $m>1$, the element
$\left(c,d\right)$ has a unique $m$-th root $\left(a,b\right)$. Indeed, let
$b=\frac{d}{m}$, and solve the equation
$$ \left(1+b+\cdots +b^{n-1}\right)a = c $$
for $a$, then clearly $\left(a,b\right)^m=\left(c,d\right)$.

\begin{prop} The one-relator group
$$ G = \gppres{x,y}{\left[\left[x,y\right],y\right]=1} $$
is not a CSA-group. However, it is $\boldq$-faithful, residually torsion-free
nilpotent, residually $p$ and with torsion-free pro-$p$-completion
for every prime $p$ (hence does not contain a nonabelian Baumslag-Solitar
group).\label{OneRelNotCSA}
\end{prop}

\underline{Proof}. Let $X$ be the normal subgroup of the free group
$F\left(x,y\right)$ generated by $x$, and let $x_i=y^{-i}xy^i$ for all
integers $i$. Then
$$ \left[\left[x,y\right],y\right] =\left[x_0^{-1}x_1,y\right]
= x_1^{-1}x_0x_1^{-1}x_2 .$$
We see that in $G$ $x_i^{-1}x_{i-1}=x_{i+1}^{-1}x_i$ for all $i\in\boldz$.
Let $d=x_1^{-1}x_0$, then $d=x_{i+1}^{-1}x_i$ for all $i\in \boldz$. It is
easy to prove, by induction on $\left|i\right|$, that $x_i=xd^{-i}$ for all
integers $i$. It follows that $X$ is freely generated by $d$ and $x$, and
$G$ is the semidirect product of the free group
$F\left(d,x\right)$ and the infinite cycle $\subgp{y}$, with $y$ acting
trivially on $d$, and $x^{y^i}=xd^{-i}$.

Let $F_k$ be the $k$-th term of the central descending series of
$F\left(d,x\right)$. Then $G/F_{k+1}$ is torsion free for all $k$ (it is a
semi-direct product of $F/F_{k+1}$ and $\subgp{y}$) We claim that it is also
nilpotent. We have $\left[y,d\right]=1$ and
$$ \left[x,d\right]^y = \left[xd^{-1},d\right] =
\left[x,d\right]^{d^{-1}}.$$
Suppose that $c$ is a commutator of weight $k\geq 2$ in $d$ and $x$. It is
``multilinear'' modulo $F_{k+1}$; hence it follows from the relations $d^y=d$
and $x^y =xd^{-1}$, that if $x$ appears only once in $c$, then
$c^y\equiv c\pmod{F_{k+1}}$ and $\left[y,c\right]\equiv 1 \pmod{F_{k+1}}$.
By the same reasoning, if $x$ appears $j\geq 2$ times in the commutator $c$,
then $\left[y,c\right]$ is congruent, modulo $F_{k+1}$, to a product of
commutators in which $x$ appears at most $j-1$ times. If we commutate the
element $c$ $j$ times by $y$, it drops into $F_{k+1}$. This shows that
$G/F_{k+1}$ is nilpotent. By Proposition~\ref{tfnilp}, $G$ is
$\boldq$-faithful.

For every prime $p$, the pro-$p$-completion of $G$
is the semi-direct product
of the free pro-$p$-group on $x,d$ by the free pro-$p$-group on the single
generator $y$. Clearly, the canonical map from $G$ into its pro-$p$-completion
is injective. Thus, $G$ is residually $p$.

That $G$ is $\boldq$-faithful follows from Proposition~\ref{tfnilp} or
Proposition~\ref{aapfaith}.

To see that $G$ is not CSA, we write $y_i=x^{-i}yx^i$, and we find that
$\left[y_0,y_1\right]=1$. Hence,
$$ 1\neq \subgp{y_1} =\subgp{y_0,y_1}\cap\subgp{y_1,y_2}
= \subgp{y_0,y_1}\cap\subgp{y_0,y_1}^x, x\notin\subgp{y_0,y_1}.$$
This completes the proof of the Proposition.

%\bibliography{proj}

\end{document}